\magnification 1200
\input plainenc
\input amssym
\fontencoding{T2A}
\inputencoding{utf-8}
\tolerance 4000
\relpenalty 10000
\binoppenalty 10000
\parindent 1.5em

\hsize 17truecm
\vsize 24truecm
\hoffset 0truecm
\voffset -0.5truecm

\font\TITLE labx1440
\font\tenrm larm1000
\font\cmtenrm cmr10
\font\tenit lati1000
\font\tenbf labx1000
\font\tentt latt1000
\font\teni cmmi10 \skewchar\teni '177
\font\tensy cmsy10 \skewchar\tensy '60
\font\tenex cmex10
\font\teneufm eufm10
\font\eightrm larm0800
\font\cmeightrm cmr8
\font\eightit lati0800
\font\eightbf labx0800
\font\eighttt latt0800
\font\eighti cmmi8 \skewchar\eighti '177
\font\eightsy cmsy8 \skewchar\eightsy '60
\font\eightex cmex8
\font\eighteufm eufm8

\font\cmsixrm cmr6

\font\sixbf labx0600
\font\sixi cmmi6 \skewchar\sixi '177
\font\sixsy cmsy6 \skewchar\sixsy '60
\font\sixeufm eufm6

\font\cmfiverm cmr5

\font\fivebf labx0500
\font\fivei cmmi5 \skewchar\fivei '177
\font\fivesy cmsy5 \skewchar\fivesy '60
\font\fiveeufm eufm5
\font\tencmmib cmmib10 \skewchar\tencmmib '177
\font\eightcmmib cmmib8 \skewchar\eightcmmib '177
\font\sevencmmib cmmib7 \skewchar\sevencmmib '177
\font\sixcmmib cmmib6 \skewchar\sixcmmib '177
\font\fivecmmib cmmib5 \skewchar\fivecmmib '177
\newfam\cmmibfam
\textfont\cmmibfam\tencmmib \scriptfont\cmmibfam\sevencmmib
\scriptscriptfont\cmmibfam\fivecmmib
\def\tenpoint{\def\rm{\fam0\tenrm}\def\it{\fam\itfam\tenit}%
	\def\bf{\fam\bffam\tenbf}\def\tt{\fam\ttfam\tentt}%
	\textfont0\cmtenrm \scriptfont0\cmsevenrm \scriptscriptfont0\cmfiverm
  	\textfont1\teni \scriptfont1\seveni \scriptscriptfont1\fivei
  	\textfont2\tensy \scriptfont2\sevensy \scriptscriptfont2\fivesy
  	\textfont3\tenex \scriptfont3\tenex \scriptscriptfont3\tenex
  	\textfont\itfam\tenit
	\textfont\bffam\tenbf \scriptfont\bffam\sevenbf
	\scriptscriptfont\bffam\fivebf
	\textfont\eufmfam\teneufm \scriptfont\eufmfam\seveneufm
	\scriptscriptfont\eufmfam\fiveeufm
	\textfont\cmmibfam\tencmmib \scriptfont\cmmibfam\sevencmmib
	\scriptscriptfont\cmmibfam\fivecmmib
	\normalbaselineskip 12pt
	\setbox\strutbox\hbox{\vrule height8.5pt depth3.5pt width0pt}%
	\normalbaselines\rm}
\def\eightpoint{\def\rm{\fam 0\eightrm}\def\it{\fam\itfam\eightit}%
	\def\bf{\fam\bffam\eightbf}\def\tt{\fam\ttfam\eighttt}%
	\textfont0\cmeightrm \scriptfont0\cmsixrm \scriptscriptfont0\cmfiverm
	\textfont1\eighti \scriptfont1\sixi \scriptscriptfont1\fivei
	\textfont2\eightsy \scriptfont2\sixsy \scriptscriptfont2\fivesy
	\textfont3\eightex \scriptfont3\eightex \scriptscriptfont3\eightex
	\textfont\itfam\eightit
	\textfont\bffam\eightbf \scriptfont\bffam\sixbf
	\scriptscriptfont\bffam\fivebf
	\textfont\eufmfam\eighteufm \scriptfont\eufmfam\sixeufm
	\scriptscriptfont\eufmfam\fiveeufm
	\textfont\cmmibfam\eightcmmib \scriptfont\cmmibfam\sixcmmib
	\scriptscriptfont\cmmibfam\fivecmmib
	\normalbaselineskip 11pt
	\abovedisplayskip 5pt
	\belowdisplayskip 5pt
	\setbox\strutbox\hbox{\vrule height7pt depth2pt width0pt}%
	\normalbaselines\rm
}

\def\No{\char 157}
\def\empty{}

\catcode`\@ 11
\catcode`\" 13
\def"#1{\ifx#1<\char 190\relax\else\ifx#1>\char 191\relax\else #1\fi\fi}

\def\newl@bel#1#2{\expandafter\def\csname l@#1\endcsname{#2}}
\openin 11\jobname .aux
\ifeof 11
	\closein 11\relax
\else
	\closein 11
	\input \jobname .aux
	\relax
\fi

\newcount\c@section
\newcount\c@subsection
\newcount\c@subsubsection
\newcount\c@equation
\newcount\c@bibl
\c@section 0
\c@subsection 0
\c@subsubsection 0
\c@equation 0
\c@bibl 0
\def\lab@l{}
\def\label#1{\immediate\write 11{\string\newl@bel{#1}{\lab@l}}%
	\ifhmode\unskip\fi}
\def\eqlabel#1{\rlap{$(\equation)$}\label{#1}}

\def\section#1{\global\advance\c@section 1
	{\par\vskip 3ex plus 0.5ex minus 0.1ex
	\rightskip 0pt plus 1fill\leftskip 0pt plus 1fill\noindent
	{\bf\S\thinspace\number\c@section .~#1}\par\penalty 25000%
	\vskip 1ex plus 0.25ex}
	\gdef\lab@l{\number\c@section.}
	\c@subsection 0
	\c@subsubsection 0
	\c@equation 0
}
\def\subsection{\global\advance\c@subsection 1
	\par\vskip 1ex plus 0.1ex minus 0.05ex{\bf\number\c@subsection. }%
	\gdef\lab@l{\number\c@section.\number\c@subsection}%
	\c@subsubsection 0\c@equation 0%
}
\def\subsubsection{\global\advance\c@subsubsection 1
	\par\vskip 1ex plus 0.1ex minus 0.05ex%
	{\bf\number\c@subsection.\number\c@subsubsection. }%
	\gdef\lab@l{\number\c@section.\number\c@subsection.%
		\number\c@subsubsection}%
}
\def\equation{\global\advance\c@equation 1
	\gdef\lab@l{\number\c@section.\number\c@subsection.%
	\number\c@equation}{\rm\number\c@equation}
}
\def\bibitem#1{\global\advance\c@bibl 1
	[\number\c@bibl]%
	\gdef\lab@l{\number\c@bibl}\label{#1}
}
\def\ref@ref#1.#2:{\def\REF@{#2}\ifx\REF@\empty{\S\thinspace#1}%
	\else\ifnum #1=\c@section {#2}\else {\S\thinspace#1.#2}\fi\fi
}
\def\ref@eqref#1.#2.#3:{\ifnum #1=\c@section\ifnum #2=\c@subsection
	{(#3)}\else{#2\thinspace(#3)}\fi\else{\S\thinspace#1.#2\thinspace(#3)}\fi
}
\def\ref#1{\expandafter\ifx\csname l@#1\endcsname\relax
	{\bf ??}\else\edef\mur@{\csname l@#1\endcsname :}%
	{\expandafter\ref@ref\mur@}\fi
}
\def\eqref#1{\expandafter\ifx\csname l@#1\endcsname\relax
	{(\bf ??)}\else\edef\mur@{\csname l@#1\endcsname :}%
	{\expandafter\ref@eqref\mur@}\fi
}
\def\cite#1{\expandafter\ifx\csname l@#1\endcsname\relax
	{\bf ??}\else\hbox{\bf\csname l@#1\endcsname}\fi
}

\def\superalign#1{\openup 1\jot\tabskip 0pt\halign to\displaywidth{%
	\tabskip 0pt plus 1fil$\displaystyle ##$&%
	\tabskip 0pt\hss $\displaystyle ##{}$&%
	$\displaystyle {}##$\hss\tabskip 0pt plus 1fil&%
	\tabskip 0pt\hss ##\crcr #1\crcr}}

\def\Wo{{\mathpalette\Wo@{}}W}
\def\Wo@#1{\setbox0\hbox{$#1 W$}\dimen@\ht0\dimen@ii\wd0\raise0.65\dimen@%
\rlap{\kern0.35\dimen@ii$#1{}^\circ$}}
\catcode`\@ 12

\def\proof{\par\medskip{\rm Д$\,$о$\,$к$\,$а$\,$з$\,$а$\,$т$\,$е$\,$л$\,$ь%
	$\,$с$\,$т$\,$в$\,$о.}\ }
\def\endproof{{\parfillskip 0pt\hfill$\square$\par}\medskip}

\immediate\openout 11\jobname.aux


\frenchspacing\rm
\leftline{УДК~517.927+517.984}\vskip 0.25cm
{\leftskip 0cm plus 1fill\rightskip 0cm plus 1fill\parindent 0cm\baselineskip 15pt
\TITLE О мажорантах собственных значений задач Штурма--Лиувилля с потенциалами
из шаров весовых пространств\par\vskip 0.15cm\rm А.$\,$А.~Владимиров%
\footnote{}{\eightpoint\rm Работа поддержана РФФИ, гранты \No\No$\,$14-01-31423
и~13-01-00705.}\par}
\vskip 0.25cm
$$
	\vbox{\hsize 0.75\hsize\leftskip 0cm\rightskip 0cm
	\eightpoint\rm
	{\bf Аннотация:\/} Изучается вопрос о точной априорной мажоранте
	наимень\-шего соб\-ственного значения задачи Штурма--Лиувилля
	$$
		\openup -1\jot
		\displaylines{-y''+qy=\lambda y,\cr y(0)=y(1)=0}
	$$
	с ограничением на потенциал вида $\int_0^1 rq^\gamma\,dx\leqslant 1$,
	где вес $r\in C(0,1)$ равномерно положителен внутри интервала $(0,1)$.
	Даётся конструктивное доказательство дости\-жимости указанной мажоранты
	для всех $\gamma>1$, а при некотором естественном расширении класса
	допустимых потенциалов~--- и для $\gamma=1$.\par
	}
$$

\vskip 0.5cm
\section{Введение}\label{par:1}%
\subsection
Зафиксируем весовую функцию $r\in C(0,1)$, равномерно положительную внутри
интервала $(0,1)$, и рассмотрим спектральную граничную задачу
$$
	\displaylines{\hbox to \displaywidth{\eqlabel{eq:1}\hfill
		$-y''+qy=\lambda y,$\hfill}\cr
		\hbox to \displaywidth{\eqlabel{eq:2}\hfill
		$y(0)=y(1)=0,$\hfill}}
$$
где потенциал принадлежит классу
$$
	A_{r,\gamma}\rightleftharpoons\left\{q\in L_{1,loc}(0,1)\;:\;
		\Bigl(q\geqslant 0\Bigr)\mathbin{\&}
		\Bigl(\int_0^1rq^\gamma\,dx\leqslant 1\Bigr)\right\},
		\leqno(\equation)
$$\label{eq:3}%
отвечающему некоторому значению $\gamma\geqslant 1$. Целью настоящей статьи является
изу\-чение вопроса о точной априорной мажоранте $M_{r,\gamma}\rightleftharpoons
\sup_{q\in A_{r,\gamma}}\lambda_0(q)$ наименьшего собственного значения
$\lambda_0(q)$ задачи \eqref{eq:1}--\eqref{eq:2}, а также о характере достижимости
этой мажоранты внутри класса $A_{r,\gamma}$ или некоторого его естественного
расширения.

\subsection
Применительно к ряду конкретных весов поставленная выше задача рассматрива\-лась
и ранее. Так, для случая $r(x)\equiv 1$ и $\gamma=1$ соответствующий результат
без доказательства сформулирован уже в работе [\cite{OEK}]. Полный разбор задачи
с весом $r(x)\equiv 1$ дан в работах [\cite{OGI}, \cite{OSZE}]. Более общий
случай веса вида $r(x)\equiv x^\alpha\cdot(1-x)^\beta$ рассматривался, в частности,
в работах [\cite{OSZT}, \cite{OSZT2}].

Подход, применённый в работах [\cite{OSZE}, \cite{OSZT}, \cite{OSZT2}],
допускает перенесение и на случай произ\-вольного веса. Он, однако, существенным
образом опирается на теорему выбора, и по\-тому неприемлем с точки зрения
конструктивного направления в математике [\cite{OKM}]. Та же теорема использована
и в работе [\cite{OGI}]. Подход последней работы, однако, допускает ряд уточнений,
позволяющих сделать доказательство полностью конструктивным. Именно это
и осуществляется в настоящей статье. Стремление к конструктивности проводимых
рассуждений обусловливает также характер приводимых доказательств ряда
вспомо\-гательных предложений (например, \ref{prop:5}).

Методологическую основу излагаемых далее результатов, как и результатов работы
[\cite{SPEE}], составляет теория операторов Штурма--Лиувилля с потенциалами-%
обобщёнными функ\-циями [\cite{OSSP}].

\subsection
Поясняющие изложение ссылки на разделы настоящей статьи далее приводятся в тексте
в прямых скобках. При этом если ссылка даётся на раздел текущего пункта (либо
параграфа), то указание на номер этого пункта (либо параграфа) опускается.

Все функциональные пространства далее предполагаются вещественными.


\section{Предварительные утверждения}\label{par:2}%
\subsection
Далее мы будем рассматривать множество обобщённых функций $\Wo_{2,loc}^{-1}(0,1)$
в качестве пространства Фреше с отвечающими всевозможным натуральным значениям
$\ell\geqslant 2$ полунормами
$$
	\|q\|_\ell\rightleftharpoons\sup_{\scriptstyle y\in\Wo_2^1(2^{-\ell},
		1-2^{-\ell})\atop\scriptstyle\|y'\|_{L_2[0,1]}\leqslant 1}
		\langle q,y\rangle.\leqno(\equation)
$$\label{eq:Frechet}%
Множество {\it неотрицательных\/} обобщённых функций класса $\Wo_{2,loc}^{-1}(0,1)$
мы далее будем обоз\-начать символом $\frak P$. С каждой обобщённой функцией $q\in
\frak P$ может быть связано гильбер\-тово пространство ${\frak D}_q$, получаемое
пополнением семейства $C_0^1(0,1)$ относительно скалярного произведения
$$
	\langle y,z\rangle_{\frak D_q}\rightleftharpoons\int_0^1 y'z'\,dx+
		\langle q,yz\rangle.\leqno(\equation)
$$\label{eq:norm}%

\subsubsection\label{prop:93}
{\it При любом выборе потенциала $q\in\frak P$, значения $\varepsilon\in (0,1/2)$
и функции $y\in\frak D_q$, ортого\-нальной подпространству $\Wo_2^1(\varepsilon,
1-\varepsilon)$, справедливо неравенство
$$
	\|y\|_{L_2[0,1]}\leqslant\sqrt\varepsilon\,\|y\|_{\frak D_q}.
$$
}%
\proof
Зафиксируем неубывающую функцию $Q\in L_{2,loc}(0,1)$ со свой\-ством
$$
	(\forall z\in C_0^1(0,1))\qquad\langle q,z\rangle=-\int_0^1 Qz'\,dx.
	\leqno(\equation)
$$\label{eq:Q}%
Независимо от выбора функции $z\in\Wo_2^1(\varepsilon,1-\varepsilon)$ выполняется
равенство
$$
	\hbox to \displaywidth{\hss $\displaystyle\int_\varepsilon^{1-\varepsilon}
		 [y'-Qy]\,z'\,dx=\int_\varepsilon^{1-\varepsilon} Qy'z\,dx$,\hss
		[\eqref{eq:norm}, \eqref{eq:Q}]}
$$
что означает подчинение функции $y$ на интервале $(\varepsilon,1-\varepsilon)$
уравнению
$$
	-[y'-Qy]'-Qy'=0.\leqno(\equation)
$$\label{eq:ekv}%
Тем самым, для любых двух точек непрерывности $a\in (\varepsilon,1-\varepsilon)$
и $b\in (a,1-\varepsilon)$ функции $Q$ выполняются соотношения
$$
	\superalign{&0&=\int_a^b [y'-Qy]'\,y\,dx+\int_a^b Qy'y\,dx
			&[\eqref{eq:ekv}]\cr
		&&=[y'-Qy]\,y\Bigr|_a^b-\int_a^b (y')^2\,dx+
			\int_a^b Q\cdot(y^2)'\,dx\cr
		&&=\left.{(y^2)'\over 2}\right|_a^b-\int_a^b (y')^2\,dx-
			\int_a^by^2\,dQ\cr
		&&\leqslant\left.{(y^2)'\over 2}\right|_a^b,\cr
	}
$$
означающие неубывание функции $(y^2)'$ на некотором почти полном подмножестве
интервала $(\varepsilon,1-\varepsilon)$. Обусловленная этим фактом вогнутость
функции $y^2$ гарантирует справедливость равенства
$$
	\sup_{x\in [\varepsilon,1-\varepsilon]} |y(x)|=
		\sup\bigl\{|y(\varepsilon)|,\,|y(1-\varepsilon)|\bigr\},
$$
с учётом выполнения вне отрезка $[\varepsilon,1-\varepsilon]$ тривиальных
неравенств $|y(x)|\leqslant\sqrt\varepsilon\,\|y'\|_{L_2[0,1]}$ автоматически
означающего искомое.
\endproof

В частности, при любом $\varepsilon\in (0,1/2)$ всякая \hbox{$\varepsilon$-сеть}
\hbox{$L_2[0,1]$-об}\-ра\-за единичного ша\-ра пространства $\Wo_2^1(\varepsilon^2,
1-\varepsilon^2)$ будет являться \hbox{$2\varepsilon$-сетью} \hbox{$L_2[0,1]$-об}%
\-ра\-за единичного ша\-ра пространства $\frak D_q$. Это замечание означает полную
непрерывность вложения $\frak D_q\hookrightarrow L_2[0,1]$.

\subsection
Пространству ${\frak D}_q$ может быть сопоставлено сопряжённое гильбертово
простран\-ство ${\frak D}_q^*$, а также линейный операторный пучок $T_q\colon\Bbb R
\to{\cal B}({\frak D}_q,{\frak D}_q^*)$ вида
$$
	\langle T_q(\lambda)y,z\rangle\rightleftharpoons\int_0^1 \bigl[y'z'-
		\lambda yz\bigr]\,dx+\langle q,yz\rangle.\leqno(\equation)
$$\label{eq:Tq}%
Упомянутая ранее полная непрерывность вложения $\frak D_q\hookrightarrow L_2[0,1]$
гарантирует, что спектр пучка $T_q$ допускает представление в виде
последовательности сосчитанных с учётом кратности собственных значений
$$
	\lambda_0(q)\leqslant\lambda_1(q)\leqslant\ldots\leqslant\lambda_n(q)
		\leqslant\ldots
$$

\subsubsection\label{prop:2}
{\it При любом выборе индекса $n\in\Bbb N$ и потенциала $q\in\frak P$ справедливо
неравенство
$$
	\lambda_n(q)\leqslant 4\pi^2\cdot(n+1)^2\cdot(1+2\,\|q\|_2).
$$
}%
\proof
Независимо от выбора функции $y\in\Wo_2^1(1/4,3/4)$ и значения $\lambda\in\Bbb R$
справедливы оценки
$$
	\superalign{&\langle T_q(\lambda)y,y\rangle&\leqslant\int_0^1\bigl[
			(y')^2-\lambda y^2\bigr]\,dx+\|q\|_2\cdot
			\sqrt{\int_0^1 4y^2\cdot(y')^2\,dx}&[\eqref{eq:Tq},
			\eqref{eq:Frechet}]\cr
		&&\leqslant\int_0^1\bigl[(1+2\,\|q\|_2)\cdot(y')^2-
			\lambda y^2\bigr]\,dx.\cr
	}
$$
Соответственно, при всяком $\lambda>4\pi^2\cdot(n+1)^2\cdot(1+2\,\|q\|_2)$
заведомо найдётся \hbox{$(n+1)$-мер}\-ное под\-пространство $\frak M\subseteq
\Wo_2^1(1/4,3/4)$ со свойством
$$
	(\forall y\in\frak M\setminus\{0\})\qquad\langle T_q(\lambda)y,y\rangle<0.
$$
Согласно известным вариационным принципам для собственных значений самосопря\-жённых
вполне непрерывных операторов [\cite{LFA}:~п.$\,$95], это как раз и означает
справедливость доказываемого предложения.
\endproof

\subsubsection\label{prop:1}
{\it При любом выборе индекса $n\in\Bbb N$ заданное на множестве $\frak P$
отображение $q\mapsto 1/\lambda_n(q)$ является равномерно непрерывным.
}%
\proof
Согласно упомянутым вариационным принципам [\cite{LFA}:~п.$\,$95], имеет место
равенство
$$
	{1\over\lambda_n(q)}=\inf_{\scriptstyle\frak M\subseteq\frak D_q\atop
		\scriptstyle\mathop{\rm codim}\frak M\leqslant n}
		\sup_{\scriptstyle y\in\frak M\atop\scriptstyle\|y\|_{\frak D_q}
		\leqslant 1}\|y\|_{L_2[0,1]}^2.
$$
Соответственно [\ref{prop:93}], при всех $\ell\geqslant 2$ справедливы оценки
$$
	{1\over\lambda_n(q)}-2^{-\ell}\leqslant\inf_{\scriptstyle\frak M\subseteq
		\Wo_2^1(2^{-\ell},1-2^{-\ell})\atop\scriptstyle\mathop{\rm codim}
		\frak M\leqslant n}\sup_{\scriptstyle y\in\frak M\atop\scriptstyle
		\|y\|_{\frak D_q}\leqslant 1}\|y\|_{L_2[0,1]}^2\leqslant
		{1\over\lambda_n(q)}.\leqno(\equation)
$$\label{eq:94}%
Далее, для любых двух потенциалов $q_1,q_2\in\frak P$ и функции
$y\in\Wo_2^1(2^{-\ell},1-2^{-\ell})$ выпол\-няются соотношения
$$
	\superalign{&\|y\|_{\frak D_{q_1}}^2&=\|y\|_{\frak D_{q_2}}^2+
			\langle q_1-q_2,y^2\rangle&[\eqref{eq:norm}]\cr
		&&\leqslant\|y\|_{\frak D_{q_2}}^2+\|q_1-q_2\|_\ell\cdot
			\|2yy'\|_{L_2[0,1]}&[\eqref{eq:Frechet}]\cr
		&&\leqslant\|y\|_{\frak D_{q_2}}^2+2\,\|q_1-q_2\|_\ell\cdot
			\|y'\|_{L_2[0,1]}^2\cr
		&&\leqslant \bigl(1+2\,\|q_1-q_2\|_\ell\bigr)\cdot
			\|y\|_{\frak D_{q_2}}^2.&[\eqref{eq:norm}]\cr
	}
$$
Отсюда и из \eqref{eq:94} вытекает справедливость оценок
$$
	{1\over\lambda_n(q_2)}-2^{-\ell}\leqslant {1+2\,\|q_1-q_2\|_\ell\over
		\lambda_n(q_1)},
$$
ввиду тривиального неравенства $\lambda_n(q_1)\geqslant\pi^2$ означающих,
что при $\|q_1-q_2\|_\ell<2^{-\ell-1}\pi^2$ расстояние между величинами
$1/\lambda_n(q_1)$ и $1/\lambda_n(q_2)$ мажорируется величиной $2^{-\ell+1}$.
\endproof

\subsection
Введём в рассмотрение функции
$$
	u_{q,n}\rightleftharpoons (y_{q,n}')^2+\lambda_n(q)y_{q,n}^2,
	\leqno(\equation)
$$\label{eq:u1}%
где через $y_{q,n}\in\frak D_q$ обозначены нормированные условием
$\|y_{q,n}\|_{L_2[0,1]}=1$ собственные функции пучка $T_q$, отвечающие собственным
значениям $\lambda_n(q)$. Имеют место следующие четыре факта.

\subsubsection\label{prop:4.1}
{\it При любом выборе индекса $n\in\Bbb N$, неотрицательного финитного потенциала
$q\in C(0,1)$ и натурального значения $\ell\geqslant 2$ справедливо соотношение
$$
	(\forall x_1,x_2\in [2^{-\ell},1-2^{-\ell}])\qquad {u_{q,n}(x_1)\over
		u_{q,n}(x_2)}\leqslant\exp\left[2^{\ell/2}\,\|q\|_{\ell+1}\right].
$$
}%
\proof
Положим $\varepsilon\rightleftharpoons 2^{-\ell-1}$. Заметим, что соотношения
$$
	\superalign{&\biggl|\bigl[(y_{q,n}')^2+\lambda_n(q)y_{q,n}^2\bigr]'\biggr|&=
			q\cdot|2y_{q,n}'y_{q,n}|&[\eqref{eq:1}]\cr
		&&\leqslant {q\over\sqrt{\lambda_n(q)}}\cdot
			\bigl[(y_{q,n}')^2+\lambda_n(q)y_{q,n}^2\bigr]\cr
	}
$$
вместе с тривиальным неравенством $\lambda_n(q)\geqslant\pi^2$ немедленно влекут
за собой оценки
$$
	(\forall x_1,x_2\in [2\varepsilon,1-2\varepsilon])\qquad{u_{q,n}(x_1)\over
		u_{q,n}(x_2)}\leqslant\exp\left[{\int_{2\varepsilon}^{1-
		2\varepsilon}q\,dx\over\pi}\right].
$$
Рассматривая функцию $w\in\Wo_2^1(\varepsilon,1-\varepsilon)$ вида
$$
	w(x)\rightleftharpoons\cases{1&при $x\in [2\varepsilon,1-2\varepsilon]$,\cr
		(x-\varepsilon)/\varepsilon&при $x\in [\varepsilon,
			2\varepsilon]$,\cr
		(1-x-\varepsilon)/\varepsilon&при $x\in [1-2\varepsilon,
			1-\varepsilon]$,}
$$
убеждаемся также в справедливости оценок
$$
	\postdisplaypenalty 10000
	\superalign{&\int_{2\varepsilon}^{1-2\varepsilon}q\,dx&\leqslant
		\int_0^1 qw\,dx\cr
		&&\leqslant\|q\|_{\ell+1}\cdot\sqrt{2\over\varepsilon}.&
		[\eqref{eq:Frechet}]\cr
	}
$$
Тем самым, доказываемое предложение является верным.
\endproof

\subsubsection\label{prop:4.2}
{\it При любом выборе индекса $n\in\Bbb N$, неотрицательного финитного потенциала
$q\in C(0,1)$ и натурального значения $\ell>4+n+2\,\|q\|_2$ справедливо соотношение
$$
	(\forall x\in [2^{-\ell},1-2^{-\ell}])\qquad u_{q,n}(x)\geqslant
		\lambda_n(q)\cdot\exp\left[-1-2^{\ell/2}\,\|q\|_{\ell+1}\right].
$$
}%
\proof
Положим $\varepsilon\rightleftharpoons 2^{-\ell}$. Характер выбора значения
$\ell$ гарантирует [\ref{prop:2}] справедливость оценки $\lambda_n(q)\varepsilon^2
<1/2$. Соответственно, выполняются соотношения
$$
	\superalign{&(\forall\zeta\in [0,1])\qquad |y_{q,n}(\zeta)|&\leqslant\left|
			(1-\zeta)\cdot\int_0^\zeta y_{q,n}'\,dx-\zeta\cdot
			\int_\zeta^1 y_{q,n}'\,dx\right|\cr
		&&\leqslant\sqrt{\zeta\cdot(1-\zeta)}\cdot\|y_{q,n}'\|_{L_2[0,1]}\cr
		\eqlabel{eq:102}&&\leqslant\sqrt{\lambda_n(q)\cdot\zeta\cdot
			(1-\zeta)},&[\eqref{eq:Tq}]\cr
		&(\forall\zeta\in [\varepsilon,1-\varepsilon])\qquad
			\lambda_n(q)&=\int_0^1\lambda_n(q)y_{q,n}^2\,dx\cr
		&&\leqslant\int_{[0,\varepsilon]\cup[1-\varepsilon,1]}
			\lambda_n(q)y_{q,n}^2\,dx+\int_\varepsilon^{1-\varepsilon}
			u_{q,n}\,dx&[\eqref{eq:u1}]\cr
		&&\leqslant[\lambda_n(q)]^2\varepsilon^2+\int_\varepsilon^{1-
			\varepsilon}u_{q,n}\,dx&[\eqref{eq:102}]\cr
		\eqlabel{eq:u0}&&\leqslant{\lambda_n(q)\over 2}+u_{q,n}(\zeta)\cdot
			\exp\left[2^{\ell/2}\,\|q\|_{\ell+1}\right],&
			[\ref{prop:4.1}]\cr
		&(\forall\zeta\in [\varepsilon,1-\varepsilon])\qquad
			u_{q,n}(\zeta)&\geqslant{\lambda_n(q)\over 2}\cdot
			\exp\left[-2^{\ell/2}\,\|q\|_{\ell+1}\right].&
			[\eqref{eq:u0}]\cr
	}
$$
Тем самым, доказываемое предложение является верным.
\endproof

\subsubsection\label{prop:4.4}
{\it При любом выборе индекса $n\in\Bbb N$, неотрицательного финитного потенциала
$q\in C(0,1)$ и натурального значения $\ell>4+n+2\,\|q\|_2$ всякая точка $x\in
[2^{-\ell+1},1-2^{-\ell+1}]$ со свойством
$$
	|y_{q,n}(x)|\leqslant 2^{-\ell}\cdot\exp\left[{-1-2^{-\ell/2}\,
		\|q\|_{\ell+1}\over 2}\right]\leqno(\equation)
$$\label{eq:u2}%
удалена от некоторого нуля функции $y_{q,n}$ не более, чем на $2^{-\ell}$.
}%
\proof
Рассмотрим произвольный отрезок $\Delta\subseteq [2^{-\ell},1-2^{-\ell}]$,
на котором функция $y_{q,n}$ является знакоопределённой, и все точки которого
обладают свойством~\eqref{eq:u2}. Тогда [\ref{prop:4.2}, \eqref{eq:u1}] функция
$y_{q,n}'$ на отрезке $\Delta$ также является знакоопре\-делённой, причём
оценивается снизу по абсолютной величине постоянной
$$
	{\sqrt{\lambda_n(q)}\over 2}\cdot\exp\left[{-1-2^{-\ell/2}\,\|q\|_{\ell+1}
		\over 2}\right].
$$
Это автоматически означает, что длина отрезка $\Delta$ не может превосходить
величину $2^{-\ell+1}/\sqrt{\lambda_n(q)}$. Учёт тривиальной оценки $\lambda_n(q)
\geqslant\pi^2$ завершает доказательство.
\endproof

\subsubsection\label{prop:4.5}
{\it При любом выборе индекса $n\in\Bbb N$, неотрицательного финитного
потенциала $q\in C(0,1)$ и натурального значения $\ell>5+n+2\,\|q\|_2$
справедливо неравенство
$$
	\postdisplaypenalty 10000
	\lambda_{n+1}(q)-\lambda_n(q)\geqslant 2^{-\ell}\cdot
		\exp\bigl[-2^{-\ell/2}\,\|q\|_{\ell+1}\bigr].
$$
}%
\proof
Пусть $a$ и $b>a$~--- два нуля функции $y_{q,n}$. Соотношения
$$
	\superalign{&0&=\int_a^b\{-y_{q,n}''+[q-\lambda_n(q)]y_{q,n}\}\cdot
			y_{q,n}\,dx&[\eqref{eq:1}]\cr
		&&\geqslant\int_a^b [(y_{q,n}')^2-\lambda_n(q)y_{q,n}^2]\,dx\cr
		&&\geqslant\left[{\pi^2\over (b-a)^2}-\lambda_n(q)\right]\cdot
			\int_a^b y_{q,n}^2\,dx\cr
	}
$$
означают, что расстояние между точками $a$ и $b$ не может быть меньше величины
$\pi/\sqrt{\lambda_n(q)}$. Аналогичным образом, расстояние между двумя различными
нулями функции $y_{q,n+1}$ не может быть меньше величины $\pi/\sqrt{%
\lambda_{n+1}(q)}$. Поскольку при этом число нулей функции $y_{q,n+1}$ превосходит
таковое для функции $y_{q,n}$, то заведомо должен найтись нуль $\zeta\in (0,1)$
функции $y_{q,n+1}$, удалённый от каждого из нулей функции $y_{q,n}$ более чем
на $\pi/[3\sqrt{\lambda_{n+1}(q)}]$.

Далее, характер выбора значения $\ell$ гарантирует справедливость
оценок
$$
	\superalign{&{\pi\over 3\sqrt{\lambda_{n+1}(q)}}&\geqslant
		{1\over 6\,(n+2)\,\sqrt{1+2\,\|q\|_2}}&[\ref{prop:2}]\cr
		&&>2^{-\ell+1},\cr
	}
$$
а потому [\ref{prop:4.4}] и оценки
$$
	|y_{q,n}(\zeta)|\geqslant 2^{-\ell}\cdot\exp\left[{-1-2^{-\ell/2}\,
		\|q\|_{\ell+1}\over 2}\right].\leqno(\equation)
$$\label{eq:u3}%
Легко также усматривается [\ref{prop:4.2}, \eqref{eq:u1}] справедливость оценки
$$
	|y_{q,n+1}'(\zeta)|\geqslant\pi\cdot\exp\left[{-1-2^{-\ell/2}\,
		\|q\|_{\ell+1}\over 2}\right].\leqno(\equation)
$$\label{eq:u4}%
Тем самым, соотношения
$$
	\superalign{&\lambda_{n+1}(q)-\lambda_n(q)&\geqslant\left|
			[\lambda_{n+1}(q)-\lambda_n(q)]\cdot\int_0^\zeta
			y_{q,n} y_{q,n+1}\,dx\right|\cr
		&&=\left|\bigl(y_{q,n}'y_{q,n+1}-y_{q,n}y_{q,n+1}'\bigr)
			\Bigr|_0^{\zeta}\right|&[\eqref{eq:1}]\cr
		&&=\left|y_{q,n}(\zeta)y'_{q,n+1}(\zeta)\right|\cr
	}
$$
означают [\eqref{eq:u3}, \eqref{eq:u4}] справедливость доказываемого предложения.
\endproof

\smallskip
Установленные в настоящем пункте оценки собственных пар $\{\lambda_n(q),y_{q,n}\}$
являются довольно грубыми. Однако для достижения основных целей настоящей статьи
они вполне достаточны. В частности, предложения \ref{prop:4.5} и \ref{prop:1}
гарантируют простоту соб\-ственных значений $\lambda_n(q)$ при любом выборе
потенциала $q\in\frak P$.

\subsection
Обозначим теперь через $\Gamma_{r,\gamma}\subseteq\frak P$ замыкание множества
$A_{r,\gamma}$ в пространстве $\Wo_{2,loc}^{-1}(0,1)$. Из предложения \ref{prop:1}
немедленно вытекает справедливость равенства
$$
	M_{r,\gamma}=\sup_{q\in\Gamma_{r,\gamma}}\lambda_0(q).
$$

\subsubsection\label{prop:5}
{\it При любом выборе значения $\gamma\geqslant 1$ множество $\Gamma_{r,\gamma}$
является компактным в простран\-стве $\Wo_{2,loc}^{-1}(0,1)$.
}%
\proof
Зафиксируем произвольным образом натуральное $\ell\geqslant 2$ и веще\-ственное
$\varepsilon>0$, а также вещественное $\delta>0$ со свойством
$$
	2\delta\cdot\bigl(1+\delta+\kern -2mm\sup_{x\in [2^{-\ell},1-2^{-\ell}]}
		r^{-1/\gamma}(x)\bigr)<\varepsilon.\leqno(\equation)
$$\label{eq:709}%
Зададим на отрезке $[2^{-\ell},1-2^{-\ell}]$ разбиение единицы $\{\chi_k\}_{k=1}^N$
индикаторами интервалов, обладающее следующими свойствами:
$$
	\belowdisplayskip 0cm
	\vcenter{\hsize 0.75\hsize\leftskip 0cm\rightskip 0cm
	При любом $k\in\overline{1,N}$ интервал с индикатором $\chi_k$ непуст,
	и его длина мажорируется величиной $\delta^2$.}
	\leqno(\equation)
$$\label{eq:711}%
$$
	\vcenter{\hsize 0.75\hsize\leftskip 0cm\rightskip 0cm
	Может быть указана линейная комбинация $g$ функций из набора
	$\{\chi_k\}_{k=1}^N$, почти всюду на отрезке $[2^{-\ell},1-2^{-\ell}]$
	приближа\-ющая функцию $r^{-1/\gamma}$ с точностью $\delta$.}
	\leqno(\equation)
$$\label{eq:712}%
Всякой функции $q\in A_{r,\gamma}$ может быть поставлена в соответствие функция
$$
	\tilde q\rightleftharpoons\sum_{k=1}^N{\int_0^1 r^{1/\gamma}q\chi_k\,dx
		\over\int_0^1\chi_k\,dx}\cdot r^{-1/\gamma}\chi_k,\leqno(\equation)
$$\label{eq:700}%
удовлетворяющая соотношениям
$$
	\belowdisplayskip 0cm
	\superalign{&\int_0^1 r\tilde q^\gamma\,dx&=\sum_{k=1}^N{\left(
			\int_0^1 r^{1/\gamma}q\chi_k\,dx\right)^\gamma\over\left(
			\int_0^1\chi_k\,dx\right)^\gamma}\cdot\int_0^1\chi_k\,dx\cr
		&&\leqslant\sum_{k=1}^N{\left(\int_0^1 rq^\gamma\chi_k\,dx\right)
			\cdot\left(\int_0^1\chi_k\,dx\right)^{\gamma-1}\over
			\left(\int_0^1\chi_k\,dx\right)^{\gamma-1}}\cr
		\eqlabel{eq:701}&&\leqslant 1,&[\eqref{eq:3}]\cr
	}
$$
$$
	\abovedisplayskip 0cm\belowdisplayskip 0cm
	\hbox to \displaywidth{\eqlabel{eq:702}\hss$\displaystyle (\forall k\in
		\overline{1,N})\qquad\int_0^1 r^{1/\gamma}\cdot(q-\tilde q)
		\chi_k\,dx=0.$\hss [\eqref{eq:700}]}
$$
Сопоставив теперь произвольной функции $y\in\Wo_2^1(2^{-\ell},1-2^{-\ell})$
линейную комбинацию $\tilde y$ функций из набора $\{\chi_k\}_{k=1}^N$, почти всюду
приближающую функцию $y$ с точностью $\delta\,\|y'\|_{L_2[0,1]}$ [\eqref{eq:711}],
устанавливаем справедливость соотношений
$$
	\superalign{&\langle q-\tilde q,y\rangle&=\int_0^1 r^{1/\gamma}\cdot
			(q-\tilde q)\cdot(r^{-1/\gamma}-g)y\,dx+
			\int_0^1 r^{1/\gamma}\cdot(q-\tilde q)\,gy\,dx\cr
		&&\leqslant 2\delta\,\|y'\|_{L_2[0,1]}+\int_0^1 r^{1/\gamma}
			\cdot(q-\tilde q)\,gy\,dx&[\eqref{eq:3}, \eqref{eq:701},
			\eqref{eq:712}]\cr
		&&=2\delta\,\|y'\|_{L_2[0,1]}+\int_0^1 r^{1/\gamma}\cdot
			(q-\tilde q)\,g\cdot (y-\tilde y)\,dx&[\eqref{eq:702}]\cr
		&&\leqslant 2\delta\,\|y'\|_{L_2[0,1]}+	2\cdot\bigl(\delta+
			\kern -2mm\sup_{x\in [2^{-\ell},1-2^{-\ell}]}
			r^{-1/\gamma}(x)\bigr)\cdot\delta\,\|y'\|_{L_2[0,1]}&
			[\eqref{eq:3}, \eqref{eq:701}, \eqref{eq:712}]\cr
		&&\leqslant\varepsilon\,\|y'\|_{L_2[0,1]}.&[\eqref{eq:709}]\cr
	}
$$
Соответственно, всякая \hbox{$\varepsilon$-сеть} пересечения множества
$A_{r,\gamma}$ с конечномерной линейной оболоч\-кой набора $\{r^{-1/\gamma}
\chi_k\}_{k=1}^N$ относительно полунормы \eqref{eq:Frechet} является
\hbox{$2\varepsilon$-сетью} мно\-жества $A_{r,\gamma}$ относительно той же полунормы.
\endproof

\subsection
Через $B_{r,\gamma,\delta}$ мы далее будем обозначать множество финитных непрерывных
потенциалов $q\in A_{r,\gamma}$, подчиняющихся оценке $\lambda_0(q)\geqslant\pi^2+
\delta$.

\subsubsection\label{prop:47}
{\it При любом выборе значения $\gamma>1$ может быть указана величина $\delta>0$
со свойством $B_{r,\gamma,\delta}\neq\varnothing$.
}%
\proof
Зафиксируем финитный непрерывный потенциал $q\in A_{r,\gamma}$, не равный нулю
как элемент пространства $L_1[0,1]$. Отвечающая ему собственная функция $y_{q,0}\in
\frak D_q$ не является собственной функцией пучка $T_0$, а потому удовлетворяет
соот\-ношениям
$$
	\superalign{&\langle T_q(\pi^2)y_{q,0},y_{q,0}\rangle&=
			\int_0^1\left[(y_{q,0}')^2+(q-\pi^2)\,y_{q,0}^2\right]\,dx&
			[\eqref{eq:Tq}]\cr
		&&\geqslant\int_0^1\left[(y_{q,0}')^2-\pi^2y_{q,0}^2\right]\,dx\cr
		&&>0.\cr
	}
$$
Отсюда немедленно вытекает справедливость неравенства $\lambda_0(q)>\pi^2$, а тогда
и суще\-ствование величины $\delta>0$ со свойством $q\in B_{r,\gamma,\delta}$.
\endproof

\subsubsection\label{prop:3}
{\it При любом выборе значений $\gamma>1$ и $\delta>0$ может быть указана
постоянная $\varkappa>0$, независимо от выбора потенциала $q\in B_{r,\gamma,
\delta}$ удовлетворяющая неравенству
$$
	\int_0^1 (q-\varkappa)\cdot y_{q,0}^2\,dx\geqslant 0.
$$
}%
\proof
Зафиксируем параметризованное значениями $\varepsilon\in (0,1/8)$ семейство
функций $w_\varepsilon\in\Wo_2^1(0,1)$ вида
$$
	w_\varepsilon(x)\rightleftharpoons\cases{0&при $x\not\in [\varepsilon,
		1-\varepsilon]$,\cr\displaystyle\cos{\pi\cdot(x-1/2)
		\over 1-2\varepsilon}&иначе.}\leqno(\equation)
$$\label{eq:99}%
Положив
$$
	\hbox to \displaywidth{\eqlabel{eq:M}\hss$\displaystyle
		\Theta\rightleftharpoons 4\pi^2\cdot\Bigl(1+2\,
		\sup_{q\in\Gamma_{r,\gamma}}\|q\|_2\Bigr)$,\hss[\ref{prop:5}]}
$$
устанавливаем при $\varepsilon\to +0$ равномерные по $q\in B_{r,\gamma,\delta}$
оценки
$$
	\belowdisplayskip 0cm
	\superalign{&\langle qw_\varepsilon,w_\varepsilon\rangle&=
		\int_\varepsilon^{1-\varepsilon} rq\cdot{w_\varepsilon^2\over r}\,
			dx\cr
		&&\leqslant\left[\int_\varepsilon^{1-\varepsilon}r\cdot
			\left({w_\varepsilon^2\over r}\right)^{\gamma/(\gamma-1)}
			dx\right]^{(\gamma-1)/\gamma}&[\eqref{eq:3}]\cr
		&&\leqslant\left[\int_\varepsilon^{1-\varepsilon}r^{-1/(\gamma-1)}\,
			dx\right]^{(\gamma-1)/\gamma}&[\eqref{eq:99}]\cr
		\eqlabel{eq:101}&&\leqslant\sup_{x\in [\varepsilon,1-\varepsilon]}
			r^{-1/\gamma}(x),\cr
	}
$$
$$
	\abovedisplayskip 0cm\belowdisplayskip 0cm
	\superalign{&\int_0^1 |y_{q,0}|\,w_\varepsilon\,dx&\geqslant
			\int_{\varepsilon^{1/3}}^{1-\varepsilon^{1/3}}
			|y_{q,0}|\,w_\varepsilon\,dx\cr
		&&\geqslant\inf_{x\in [\varepsilon^{1/3},1-\varepsilon^{1/3}]}
			{w_\varepsilon(x)\over |y_{q,0}(x)|}\cdot
			\int_{\varepsilon^{1/3}}^{1-\varepsilon^{1/3}}
			y_{q,0}^2\,dx\cr
		&&\geqslant{\pi\cdot(\varepsilon^{1/3}-\varepsilon)\over
			(1-2\varepsilon)\cdot\sqrt\Theta}\cdot
			[1-\Theta\varepsilon^{2/3}]&[\ref{prop:2}, \eqref{eq:99},
			\eqref{eq:102}, \eqref{eq:M}]\cr
		\eqlabel{eq:103}&&\geqslant{\pi\varepsilon^{1/3}\over\sqrt\Theta}
			\cdot [1+o(1)],\cr
	}
$$
$$
	\abovedisplayskip 0cm
	\superalign{&\sqrt{\langle qy_{q,0},y_{q,0}\rangle}&\geqslant
			{\langle q\,|y_{q,0}|,w_\varepsilon\rangle\over
			\sqrt{\langle qw_\varepsilon,
			w_\varepsilon\rangle}}\cr
		&&={\langle |y_{q,0}|''+\lambda_0(q)|y_{q,0}|,w_\varepsilon\rangle
			\over\sqrt{\langle qw_\varepsilon,w_\varepsilon\rangle}}
			&[\eqref{eq:1}]\cr
		&&={\int_\varepsilon^{1-\varepsilon} |y_{q,0}|\cdot
			\left[w_\varepsilon''+\lambda_0(q)w_\varepsilon\right]\,dx-
			\left(|y_{q,0}|\,w_\varepsilon'\right)\Bigr|_{\varepsilon+0}^{1-
			\varepsilon-0}\over\sqrt{\langle qw_\varepsilon,
			w_\varepsilon\rangle}}\cr
		&&\geqslant{[\delta+o(1)]\cdot\int_0^1 |y_{q,0}|\,w_\varepsilon\,dx+
			O(\varepsilon^{1/2})\over\sqrt{\langle qw_\varepsilon,
			w_\varepsilon\rangle}}&[\eqref{eq:99}, \eqref{eq:102}]\cr
		&&\geqslant{\pi\delta\varepsilon^{1/3}\over\sqrt\Theta}\cdot
			\inf_{x\in [\varepsilon,1-\varepsilon]} r^{1/2\gamma}(x)
			\cdot [1+o(1)].&[\eqref{eq:103}, \eqref{eq:101}]\cr
	}
$$
Из последней оценки фиксацией достаточно малого значения $\varepsilon>0$
получается иско\-мое.
\endproof

\subsubsection\label{prop:4}
{\it При любом выборе значений $\gamma>1$, $\delta>0$ и функции $q\in B_{r,\gamma,
\delta}$ справедлива оценка
$$
	M_{r,\gamma}\geqslant\lambda_0(q)+\varkappa\cdot
		\left(1-\left[\int_0^1 rq^\gamma\,dx\right]^{1/\gamma}\right),
$$
где коэффициент $\varkappa>0$ определён предложением~\ref{prop:3}.
}%
\proof
Зафиксируем произвольное $t\geqslant 1$ со свойством $tq\in B_{r,\gamma,\delta}$.
Соотношения
$$
	\superalign{&0&\leqslant\langle T_q(\lambda_0(q))y_{tq,0},y_{tq,0}\rangle\cr
		&&=\langle T_{tq}(\lambda_0(tq))y_{tq,0},y_{tq,0}\rangle+
		\int_0^1\left[(t^{-1}-1)\cdot tq+\bigl(\lambda_0(tq)-
		\lambda_0(q)\bigr)\right]\cdot y_{tq,0}^2\,dx&[\eqref{eq:Tq}]\cr
		&&\leqslant\varkappa\cdot (t^{-1}-1)+(M_{r,\gamma}-\lambda_0(q))
		&[\ref{prop:3}]\cr
	}
$$
означают справедливость оценки $M_{r,\gamma}\geqslant\lambda_0(q)+\varkappa
\cdot(1-t^{-1})$. Последняя, в свою очередь, немед\-ленно влечёт искомое.
\endproof


\section{Основные результаты}\label{par:3}%
\subsection
Рассмотрим сначала случай $\gamma>1$. Здесь, согласно \ref{prop:5}, \ref{prop:1},
\ref{prop:47} и \ref{prop:2}, при некотором $\delta>0$ может быть указана
последовательность $\{q_n\}_{n=0}^\infty$ потенциалов класса $B_{r,\gamma,\delta}$
со свойством $\lambda_0(q_n)\to M_{r,\gamma}$. При этом, очевидно, можно
без ограничения общности считать выполненным тождество
$$
	(\forall n\in\Bbb N)\qquad\int_0^1 rq_n^\gamma\,dx=1.
$$
Рассмотрим теперь двойную последовательность $\{q_{n,m}\}_{(n,m)\in\Bbb N^2}$
потенциалов вида
$$
	q_{n,m}\rightleftharpoons {q_n+q_m\over 2}.\leqno(\equation)
$$\label{eq:qnm}%
Поскольку отображение $\lambda_0\colon L_1[0,1]\to\Bbb R$, согласно вариационным
принципам, является выпуклым, то все потенциалы $q_{n,m}$ принадлежат классу
$B_{r,\gamma,\delta}$, причём нижний предел числовой двойной последовательности
$\{\lambda_0(q_{n,m})\}_{(n,m)\in\Bbb N^2}$ не может оказаться меньшим, нежели
$M_{r,\gamma}$. Это автоматически означает справедливость равенства
$$
	\lim_{\scriptstyle n\to\infty\atop\scriptstyle  m\to\infty}
		\lambda_0(q_{n,m})=M_{r,\gamma},\leqno(\equation)
$$\label{eq:Mrg}%
а тогда, согласно предложению \ref{prop:4}, и равенства
$$
	\lim_{\scriptstyle n\to\infty\atop\scriptstyle  m\to\infty}
		\int_0^1 rq_{n,m}^\gamma\,dx=1.
$$
Известный [\cite{LFA}:~п.$\,$145] факт равномерной выпуклости весового пространства
$L_\gamma(r;\,[0,1])$ означает потому фундаментальность последовательности
$\{q_n\}_{n=0}^\infty$ в этом простран\-стве. Соответственно [\ref{prop:5},
\ref{prop:1}], справедливо следующее предложение:

\subsubsection\label{prop:3.1}
{\it При всяком $\gamma>1$ существует и однозначно определён потенциал $\hat q\in
A_{r,\gamma}$, удовлетво\-ряющий равенству $\lambda_0(\hat q)=M_{r,\gamma}$.
}%

\medskip
Зафиксируем теперь произвольным образом неотрицательную финитную функцию
$p\in C(0,1)$, а также величину
$$
	\alpha>\int_0^1 rp\hat q^{\gamma-1}\,dx-1.\leqno(\equation)
$$\label{eq:alpha}%
При любом $\varepsilon\in [0,1)$ отвечающее потенциалу
$$
	q_\varepsilon\rightleftharpoons{(1-\varepsilon)\hat q+\varepsilon p\over
		1+\alpha\varepsilon}
$$
пространство $\frak D_{q_\varepsilon}$ совпадает, с точностью до выбора
эквивалентной нормы, с простран\-ством $\frak D_{\hat q}$. Применим теперь
к параметрическому семейству пучков $T_{q_\varepsilon}\colon\Bbb R\to
{\cal B}(\frak D_{\hat q},\frak D_{\hat q}^*)$ известные результаты
[\cite{LFA}:~п.$\,$136] о поведении простых собственных значений при аналитических
возмуще\-ниях. Согласно последним, при $\varepsilon\to+0$ должна быть справед\-лива
асимптотика
$$
	\lambda_0(q_\varepsilon)=M_{r,\gamma}+\varepsilon\cdot\int_0^1 [p-(\alpha+1)
		\hat q]\cdot y_{\hat q}^2\,dx+o(\varepsilon),
$$
где через $y_{\hat q}\in\frak D_{\hat q}$ обозначена нормированная условием
$\|y_{\hat q}\|_{L_2[0,1]}=1$ собственная функция пучка $T_{\hat q}$, отвечающая
собственному значению $M_{r,\gamma}$. Поскольку при всех достаточно малых значениях
$\varepsilon>0$ справедливо соотношение $q_\varepsilon\in A_{r,\gamma}$
[\eqref{eq:alpha}, \eqref{eq:3}], это, ввиду произвольности выбора параметра
$\alpha$ со свойством \eqref{eq:alpha}, приводит к оценке
$$
	\int_0^1 rp\cdot\left[{y_{\hat q}^2\over r}-\left(\int_0^1\hat q
		y_{\hat q}^2\,dx\right)\cdot\hat q^{\gamma-1}\right]\,dx\leqslant 0.
$$
Последняя, ввиду произвольности выбора функции $p$, означает справедливость
нера\-венств
$$
	{y_{\hat q}^2(x)\over r(x)}\leqslant\left(\int_0^1\hat qy_{\hat q}^2\,dx
		\right)\cdot\hat q^{\gamma-1}(x)\qquad\hbox{при почти всех
		$x\in (0,1)$}.\leqno(\equation)
$$\label{eq:200}%
С другой стороны, всякий потенциал $\hat q\in A_{r,\gamma}$ со свойством
\eqref{eq:200} удовлетворяет соот\-ношениям
$$
	\superalign{&M_{r,\gamma}-\lambda_0(\hat q)&\leqslant\sup_{q\in
			B_{r,\gamma,0}}\left(\int_0^1[\lambda_0(q)-
			\lambda_0(\hat q)]\,y_{\hat q}^2\,dx+\langle
			T_q(\lambda_0(q))y_{\hat q},y_{\hat q}\rangle\right)\cr
		&&=\sup_{q\in B_{r,\gamma,0}}\int_0^1 [(y_{\hat q}')^2+
			(q-\lambda_0(\hat q))\,y_{\hat q}^2]\,dx
			&\llap{[\eqref{eq:Tq}]}\cr
		&&\leqslant\int_0^1 [(y_{\hat q}')^2-\lambda_0(\hat q)y_{\hat q}^2]
			\,dx+\left[\int_0^1 r\cdot\left({y_{\hat q}^2\over r}
			\right)^{\gamma/(\gamma-1)}\,dx\right]^{(\gamma-1)/\gamma}
			&\llap{[\eqref{eq:3}]}\cr
		\eqlabel{eq:202}&&\leqslant\int_0^1 [(y_{\hat q}')^2-
			\lambda_0(\hat q)y_{\hat q}^2]\,dx+
			\int_0^1\hat qy_{\hat q}^2\,dx&\llap{[\eqref{eq:200},
			\eqref{eq:3}]}\cr
		&&=\langle T_{\hat q}(\lambda_0(\hat q))y_{\hat q},
			y_{\hat q}\rangle&\llap{[\eqref{eq:Tq}]}\cr
		\eqlabel{eq:201}&&=0,\cr
	}
$$
что означает выполнение равенства $\lambda_0(\hat q)=M_{r,\gamma}$. При этом
все неравенства внутри последней выкладки вырождаются в равенства, а потому
[\eqref{eq:202}, \eqref{eq:3}] равенствами оказываются также и неравенства
\eqref{eq:200}. Соответственно, справедливо следующее пред\-ложение:

\subsubsection\label{prop:3.1.2}
{\it Потенциал $\hat q\in A_{r,\gamma}$ со свойством $\lambda_0(\hat q)=
M_{r,\gamma}$ характеризуется условием
$$
	\hat q={r^{1/(1-\gamma)}y_{\hat q}^{2/(\gamma-1)}\over\left(\int_0^1
		r^{1/(1-\gamma)}y_{\hat q}^{2\gamma/(\gamma-1)}\,dx
		\right)^{1/\gamma}}.
$$
}%

\subsection
Перейдём теперь к рассмотрению случая $\gamma=1$. Как и в предыдущем случае,
предложения \ref{prop:5}, \ref{prop:1} и \ref{prop:2} гарантируют существование
последовательности $\{q_n\}_{n=0}^\infty$ потенциалов класса $B_{r,1,0}$,
обладающей свойством $\lambda_0(q_n)\to M_{r,1}$. Обозначим через $y_n$
нормированные условием $\|y_n\|_{L_2[0,1]}=1$ неотрицательные собственные функции
пучков $T_{q_n}$. Заметим, что при любых значениях $n,m\in\Bbb N$ нормированная
условием $\|y_{n,m}\|_{L_2[0,1]}=1$ неотрицательная собственная функция пучка
$T_{q_{n,m}}$, отвечающего потенциалу \eqref{eq:qnm}, удовлетворяет хотя бы одному
из неравенств
$$
	\leqalignno{\|y_{n,m}-y_n\|_{L_2[0,1]}&>{\|y_n-y_m\|_{L_2[0,1]}\over 3},
			&\eqlabel{eq:nmn}\cr
		\|y_{n,m}-y_m\|_{L_2[0,1]}&>{\|y_n-y_m\|_{L_2[0,1]}\over 3}.
	}
$$
Без ограничения общности можно предположить выполнение первого из них. При этом
с очевид\-ностью справедливо также неравенство
$$
	\|y_{n,m}+y_n\|_{L_2[0,1]}\geqslant\|y_{n,m}-y_n\|_{L_2[0,1]}.
	\leqno(\equation)
$$\label{eq:234}%
Зафиксировав теперь не зависящую от выбора значений $n,m\in\Bbb N$ величину $\eta>0$
достаточно малой и раскладывая функцию $y_{n,m}$ в пространстве $\frak D_{q_n}$
по собственным функциям пучка $T_{q_n}$, убеждаемся в справедливости соотношений
$$
	\belowdisplayskip 0cm
	\superalign{&\int_0^1 [(y_{n,m}')^2+q_ny_{n,m}^2]\,dx&\geqslant
			\lambda_0(q_n)\cdot\left(\int_0^1 y_{n,m}y_n\,dx
			\right)^2+\cr &&\hskip 2cm{}+\lambda_1(q_n)\cdot
			\left[1-\left(\int_0^1 y_{n,m}y_n\,dx\right)^2\right]\cr
		\eqlabel{eq:145}&&\geqslant\lambda_0(q_n)+\eta\,
			\|y_n-y_m\|_{L_2[0,1]}^2,&\llap{[\eqref{eq:nmn},
			\eqref{eq:234}, \ref{prop:5}, \ref{prop:4.5}]}\cr
	}
$$
$$
	\abovedisplayskip 0cm
	\superalign{&\lambda_0(q_{n,m})&=\int_0^1 [(y_{n,m}')^2+
			q_{n,m}y_{n,m}^2]\,dx\cr
		&&={1\over 2}\cdot\biggl(\int_0^1 [(y_{n,m}')^2+q_ny_{n,m}^2]\,dx+
			\int_0^1 [(y_{n,m}')^2+q_my_{n,m}^2]\,dx\biggr)
			&\llap{[\eqref{eq:qnm}]}\cr
		&&\geqslant{\lambda_0(q_n)+\lambda_0(q_m)\over 2}+\eta\,
			\|y_n-y_m\|_{L_2[0,1]}^2.&[\eqref{eq:145}]\cr
	}
$$
С учётом справедливости, как и для случая $\gamma>1$, равенства \eqref{eq:Mrg},
последние оценки означают фундаментальность функциональной последовательности
$\{y_n\}_{n=0}^\infty$ в простран\-стве $L_2[0,1]$. Следовательно, на плотном
подмножестве $C_0^2(0,1)\subset\Wo_2^1(0,1)$ выполняются соотношения
$$
	\eqalign{(\forall\varphi\in C_0^2(0,1))\quad
		\lim_{\scriptstyle n\to\infty\atop\scriptstyle m\to\infty}
		\int_0^1 (y'_n-y'_m)\cdot\varphi'\,dx
		&=\lim_{\scriptstyle n\to\infty\atop\scriptstyle m\to\infty}
		\int_0^1 (y_m-y_n)\cdot\varphi''\,dx\cr &=0,
	}
$$
ввиду справедливости оценок $\|y_n'\|_{L_2[0,1]}\leqslant\sqrt{M_{r,1}}$
гарантирующие слабую фундамен\-тальность последовательности $\{y_n\}_{n=0}^\infty$
в пространстве $\Wo_2^1(0,1)$. Полная непрерывность вложения $\Wo_2^1(0,1)
\hookrightarrow C[0,1]$ означает при этом сильную фундаментальность
после\-довательности $\{y_n\}_{n=0}^\infty$ в пространстве $C[0,1]$. Непрерывность
линейной биекции $y\mapsto y''$ пространства $\Wo_2^1(0,1)$ на прост\-ранство
$\Wo_2^{-1}(0,1)$ означает также слабую фун\-даментальность последовательности
функций $q_ny_n=y_n''+\lambda_0(q_n)y_n$ в пространстве $\Wo_2^{-1}(0,1)$.

Зафиксируем теперь произвольное натуральное $\ell\geqslant 2$ и положим
$\varepsilon\rightleftharpoons 2^{-\ell}$. Выпол\-нение [\eqref{eq:3}] при любом
выборе индекса $n\in\Bbb N$ и функции $z\in\Wo_2^1(\varepsilon,1-\varepsilon)$
неравенства
$$
	\left|\int_0^1 q_ny_nz\,dx\right|\leqslant L\cdot\|z\|_{C[0,1]},
	\leqno(\equation)
$$\label{eq:119}%
где положено
$$
	L>\sqrt{M_{r,1}}\cdot\sup_{x\in [\varepsilon,1-\varepsilon]} r^{-1}(x),
$$
вместе с вполне непрерывным характером вложения $\Wo_2^1(\varepsilon,1-\varepsilon)
\hookrightarrow C[0,1]$ означают фундамен\-тальность последовательности
$\{q_ny_n\}_{n=0}^\infty$ относительно полунормы \eqref{eq:Frechet}. Зафиксировав
теперь не зависящую от выбора значений $n,m\in\Bbb N$ величину $\theta>0$ достаточно
малой, убеждаемся в выполне\-нии при любом выборе функции
$z\in\Wo_2^1(\varepsilon,1-\varepsilon)$ соотношений
$$
	\superalign{&\langle q_n-q_m,z\rangle&=\int_0^1[q_ny_n-q_my_m]\cdot
			{z\over y_m}\,dx+\int_0^1 q_ny_n\cdot{y_m-y_n\over y_ny_m}
			\cdot z\,dx\cr
		&&\leqslant\|q_ny_n-q_my_m\|_\ell\cdot\left\|{z'y_m-zy_m'\over
			y_m^2}\right\|_{L_2[0,1]}+{}\cr &&\hskip 2cm {}+L\cdot
			\|z'\|_{L_2[0,1]}\cdot\sup_{x\in [\varepsilon,
			1-\varepsilon]} {y_m(x)-y_n(x)\over y_n(x)y_m(x)}&
			[\eqref{eq:Frechet}, \eqref{eq:119}]\cr
		&&\leqslant\biggl(\|q_ny_n-q_my_m\|_\ell\cdot{\theta+\sqrt{M_{r,1}}
			\over\theta^2}+{}\cr &&\hskip 2cm{}+L\theta^{-2}\cdot
			\|y_n-y_m\|_{C[0,1]}\biggr)\cdot\|z'\|_{L_2[0,1]}.
			&[\ref{prop:5}, \ref{prop:4.4}]\cr
	}
$$
Тем самым, последовательность $\{q_n\}_{n=0}^\infty$ фундаментальна относительно
полунормы \eqref{eq:Frechet}. Соответственно [\ref{prop:1}], справедливо следующее
предложение:

\subsubsection\label{prop:3.2}
{\it Существует и однозначно определён потенциал $\hat q\in\Gamma_{r,1}$,
удовлетворяющий равенству $\lambda_0(\hat q)=M_{r,1}$.
}%

\medskip
Зафиксируем теперь произвольным образом функцию $p\in B_{r,1,0}$ и рассмотрим
параметризованное значениями $\varepsilon\in [0,1)$ семейство пучков
$T_{q_\varepsilon}\colon\Bbb R\to{\cal B}(\frak D_{\hat q},\frak D_{\hat q}^*)$,
отвеча\-ющих потенциалам $q_\varepsilon\rightleftharpoons(1-\varepsilon)\,\hat q+
\varepsilon p$. Для него при $\varepsilon\to+0$ справедлива асимптотика
$$
	\lambda_0(q_\varepsilon)=M_{r,1}+\varepsilon\cdot\langle p-\hat q,
		y_{\hat q}^2\rangle+o(\varepsilon),
$$
где через $y_{\hat q}\in\frak D_{\hat q}$ обозначена нормированная условием
$\|y_{\hat q}\|_{L_2[0,1]}=1$ собственная функция пучка $T_{\hat q}$, отвечающая
собственному значению $M_{r,1}$. Это автоматически означает справедливость
неравенства
$$
	\int_0^1 rp\cdot{y_{\hat q}^2\over r}\,dx\leqslant\langle\hat q,
		y_{\hat q}^2\rangle,
$$
а тогда, ввиду произвольности выбора функции $p$, и неравенства
$$
	\sup_{x\in (0,1)}{y_{\hat q}^2(x)\over r(x)}\leqslant
		\langle\hat q,y_{\hat q}^2\rangle.\leqno(\equation)
$$\label{eq:delta}%
Повторяя теперь рассуждения, проведённые ранее при получении оценок \eqref{eq:201},
устанавливаем, что всякий потенциал $\hat q\in\Gamma_{r,1}$ со свойством
\eqref{eq:delta} удовлетворяет соот\-ношению $\lambda_0(\hat q)=M_{r,1}$, причём
само неравенство \eqref{eq:delta} вырождается в равенство. Соответственно,
справедливо следующее предложение:

\subsubsection
{\it Потенциал $\hat q\in\Gamma_{r,1}$ со свойством $\lambda_0(\hat q)=M_{r,1}$
характеризуется условием
$$
	\sup_{x\in (0,1)} {y_{\hat q}^2(x)\over r(x)}=
		\langle\hat q,y_{\hat q}^2\rangle.\leqno(\equation)
$$\label{eq:303}
}%

\medskip
Соотношения $\hat q\in\Gamma_{r,1}$ и~\eqref{eq:303} означают, что функция $r^{-1}
y_{\hat q}^2$ принимает своё максималь\-ное значение почти всюду по мере $\hat q$.
Тем самым, имеет место следующий факт:

\subsubsection
{\it Никакая точка $\zeta\in (0,1)$ со свойством
$$
	{y_{\hat q}^2(\zeta)\over r(\zeta)}<\sup_{x\in (0,1)}{y_{\hat q}^2(x)
		\over r(x)}
$$
не принадлежит носителю меры $\hat q$.
}%

\subsection
Отметим, что доказательство предложения \ref{prop:3.2} практически без изменений
переносится на случай, когда вместо класса $A_{r,1}$ рассматривается его выпуклое
подмножество (ср., например, \hbox{[\cite{OEK}: Теорема~1]}). Обсуждение возможных
дальнейших обобщений полученных результатов мы оставляем за рамками настоящей
статьи.


\vskip 0.4cm
\eightpoint\rm
{\leftskip 0cm\rightskip 0cm plus 1fill\parindent 0cm
\bf Литература\par\penalty 20000}\vskip 0.4cm\penalty 20000
\bibitem{OEK} {\it Ю.$\,$В.~Егоров, В.$\,$А.~Кондратьев.\/} Об оценках первого
собственного значения задачи Штурма--Лиувилля~// Успехи матем.~наук.~--- 1984.~---
Т.~39, \No~2.~--- С.~151--152.

\bibitem{OGI} {\it В.$\,$А.~Винокуров, В.$\,$А.~Садовничий\/}. О границах изменения
собственного значения при~изменении потен\-циала~// Доклады Акад. Наук.~--- 2003.~---
Т.~392, \No~5.~--- С.~592--597.

\bibitem{OSZE} {\it С.$\,$С.~Ежак\/}. Оценки первого собственного значения задачи
Штурма--Лиувилля с условиями Дирихле~/ В~кн.: Качественные свойства решений
дифференциальных уравнений и смежные вопросы спектрального анализа. М.:~ЮНИТИ-ДАНА,
2012. С.~517--559.

\bibitem{OSZT} {\it М.$\,$Ю.~Тельнова\/}. Оценки первого собственного значения
задачи Штурма--Лиувилля с~условиями Дирихле и весовым интегральным условием~/
В~кн.: Качественные свойства решений дифференциаль\-ных уравнений и смежные
вопросы спектрального анализа. М.:~ЮНИТИ-ДАНА, 2012. С.~608--647.

\bibitem{OSZT2} {\it М.$\,$Ю.~Тельнова\/}. Об одной оценке сверху первого
собственного значения задачи Штурма--Лиу\-вилля с~условиями Дирихле и весовым
интегральным условием~// Международная мини\-кон\-ферен\-ция "<Качественная теория
дифференциальных уравнений и приложения"> (22~июня и 19~декабря 2013~г.,
24~мая 2014~г.). Сборник трудов. М.: 2014. С.~126--140.

\bibitem{OKM} {\it А.$\,$А.~Марков\/}. О конструктивной математике~// Труды МИАН
им.~В.$\,$А.~Стеклова.~--- 1962.~--- Т.~67.~--- С.~8--14.

\bibitem{SPEE} {\it E.$\,$S.~Karulina, A.$\,$A.~Vladimirov\/}. The Sturm--Liouville
problem with singular potential and~the~extrema of the first eigenvalue~// Tatra
Mountains Math.~Publications.~--- 2013.~--- V.~54.~--- P.~101--118.

\bibitem{OSSP} {\it М.$\,$И.~Нейман-заде, А.$\,$А.~Шкаликов\/}. Операторы Шрёдингера
с сингулярными потенциалами из прост\-ранств мультипликаторов~// Матем. заметки.~---
1999.~--- Т.~66, \No~5.~--- С.~723--733.

\bibitem{LFA} {\it Ф.~Рисс, Б.~Сёкефальви--Надь\/}. Лекции по функциональному
анализу, изд.~2. М.: Мир, 1979.
\bye